\documentclass[a4paper,11pt]{amsart}       % onecolumn (second format)
\usepackage{amsmath, amssymb, pinlabel}
\usepackage{epsfig, color, graphicx, latexsym%, theorem
}
\usepackage{mathptmx}      % use Times fonts if available on your TeX system
\theoremstyle{plain}
\newtheorem{theorem}{Theorem}%[section] % 1st argument is your name for it
\newtheorem{lemma}{Lemma}     % 2nd argument is what is printed

\theoremstyle{definition}
\newtheorem{conjecture}{Conjecture}  % 2nd argument is what is printed

\newtheorem{remark}{Remark}

\definecolor{mygray}{gray}{0.5}

\newcommand{\Li}{{\rm Li}_2}
\newcommand{\boldeta}{\boldsymbol{\eta}}

\allowdisplaybreaks

\begin{document}
\title{Asymptotics of  quantum $6j$ symbols}
\author{Qingtao Chen}
\author{Jun Murakami}
%
%\affil{Division of Science,
%%\affil{Department of Mathematics,
%Faculty of Science and Engineering,
%Waseda University}
%\institute{Division of Science,
%   New York University Abu Dhabi, PO Box 129188,
%   Abu Dhabi, United Arab Emirates\\
%   \email{chenqtao@nyu.edu}
%%
%%
%\and 
%%%\givenname{Jun}
%%%\surname{Murakami}
%Department of Mathematics,
%Faculty of Science and Engineering,
%Waseda University,
%3-4-1 Ohkubo, Shinjuku-ku, Tokyo 169-8555,
%Japan\\
%\email{murakami@waseda.jp}
%}

%  First author
%
%\author{Qingtao Chen and Jun Murakami}
%\givenname{Qingtao}
%\surname{Chen}

\thanks{The first author was supported by the National Center of Competence in Research SwissMAP of the Swiss National Science Foundation, and
the second author was supported in part by KAKENHI (25287014, 16H06336, 17K18728).
}

\begin{abstract}
The asymptotics of the quantum $6j$ symbols corresponding to a hyperbolic tetrahedron is investigated and 
the first two leading terms are determined for the case that the  tetrahedron has an ideal or an ultra-ideal vertex.  
These terms are given by the volume and the determinant of the Gram matrix of the  tetrahedron.  
%\keywords{quantum $6j$ symbol \and hyperbolic tetrahedron}
%\subclass{57M27 \and 20G42}
\end{abstract}

\maketitle
%\keyword{quantum $6j$ symbol}
%\keyword{hyperbolic tetrahedron}
%\subject{primary}{msc2010}{57M27}
%\subject{secondary}{msc2010}{20G42}

%\arxivreference{arXiv:1706.04887}  %%% please supply if 
%\arxivpassword{zwfz4}   %%% paper is in the arXiv

\section{Introduction}

The quantum $6j$ symbol was introduced by A. Kirillov and N. Reshetikhin  from the representation theory of the quantum group ${\mathcal U}_q(sl_2)$ and was used to reconstruct quantum $sl_2$ invariants of knots and knotted graphs in \cite{KR}.  
Then V. Turaev and O. Viro constructed a 3-manifold invariant from its tetrahedral decomposition by assigning the quantum $6j$ symbol to each tetrahedron.  
The quantum $6j$ symbol is a quantum analogue of the classical $6j$ symbol introduced by Wigner, which is  a highly symmetric version of Racah's W-coefficients expressing the recoupling coefficient of 3 angular momenta.   
The relation between the asymtotics of the classical $6j$ symbols and the volume of an Eudlidean tetrahedron was conjectured by G. Ponzano and T. Regge \cite{PR}, and was proved by J. Roberts \cite{R1}.  
\par
\begin{figure}[htb]
\[
\begin{matrix}
\qquad
\epsfig{file=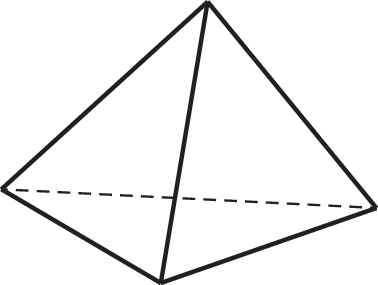, scale=0.55}
\hspace{-3.5cm}
\raisebox{11mm}{
$\begin{matrix}
&a & &&   e &\\[14pt]
 && &b  & \\[-4pt]
 & \textcolor{mygray}{c}\ \  &&& \\[10pt]
f \ \ \ && &&& \!\!\!d 
\end{matrix}$
}
\qquad&\qquad\qquad\qquad
\epsfig{file=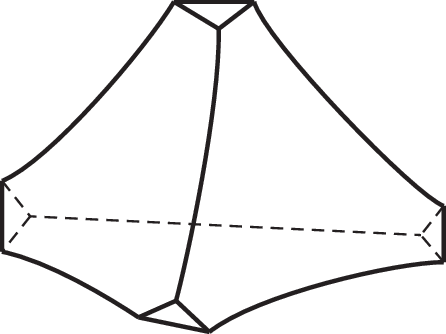, scale=0.53}
\hspace{-3.6cm}
\raisebox{11mm}{%\small
$\begin{matrix}
a& & &&  e \\[12pt]
 && b& \ \   & \\[-4pt]
 & \,\,\textcolor{mygray}{c}\ \ &&& \\[17pt]
f \ \ \ \ && &&d 
\end{matrix}\qquad$
}
\\[5pt]\qquad
\text{Usual (compact) tetrahedron}
&\qquad\qquad\qquad
\text{\quad Truncated tetrahedron}
\end{matrix}
\]
\caption{Edges of a tetrahedron and a truncated tetrahedron}
\label{fig:tet}
\end{figure}

Here we give the asymptotics of the quantum $6j$ symbols corresponding to a hyperbolic tetrahedron with at least one ideal or ultra-ideal vertex.  
Let $T$ be a hyperbolic tetrahedron, $v$ be one of its vertex, $a$, $b$, $c$ be the edges of $T$ touching to $v$ and $\theta_1$, $\theta_b$, $\theta_c$ be the dihedral angles at $a$, $b$, $c$.  
Then $\theta_a+\theta_b+\theta_c > \pi$ if $v$ is a usual vertex, $\theta_a+\theta_b+\theta_c = \pi$ if $v$ is an ideal vertex and $\theta_a+\theta_b+\theta_c < \pi$ if $v$ is a ultra-ideal vertex.  
Let $r \geq 3$ be a positive odd integer and
$
\xi_r = \exp(2\pi\sqrt{-1}/r) 
$.
A triplet $(a, b, c)$ of half integers is called {\it admissible} if it satisfies the Clebsch-Gordan condition, i.e.  $|a-b| \leq c \leq a+b$ and $a+b+c \in {\mathbf Z}$.  
\begin{theorem}
Let $T$ be a hyperbolic tetrahedron at least one of whose vertices are  ideal or ultra-ideal.  
Let  $\theta_a$, $\theta_b$, $\theta_c$, $\theta_d$, $\theta_e$, $\theta_f$ be dihedral angles of $T$ where $a$, $\cdots$, $f$ are edges of $T$ given as in Figure \ref{fig:tet}
and
let
$a_r$, $b_r$, $\cdots$, $f_r$ be sequences of half integers for odd $r\geq3$ satisfying 
\begin{equation}
\lim_{r\to\infty} \frac{a_r}{r} = \frac{1}{4}-\frac{\theta_a}{4\pi},\  
\lim_{r\to\infty} \frac{b_r}{r} = \frac{1}{4}-\frac{\theta_b}{4\pi},\   
\cdots, \ 
\lim_{r\to\infty} \frac{f_r}{r} = \frac{1}{4}-\frac{\theta_f}{4\pi}, 
\label{eq:sequence}
\end{equation}
where the triplets $(a_r, b_r, e_r)$, $(a_r, c_r, f_r)$, $(b_r, d_r, f_r)$ and $(c_r, d_r, e_r)$ are all admissible.  
We also assume that 
%$ 0 < a_r$, $b_r$, $\cdots$, $f_r < \frac{r}{4}$,  and 
$a_r+b_r+e_r \geq \frac{r-3}{2}$ for the three edges $a$, $b$, $e$  attached to the same ideal or ultra-ideal vertex.  
Then
\begin{equation}
\lim_{r\to\infty}\frac{2\, \pi}{r} 
\log\left|
\left\{\begin{matrix}
a_r & b_r & e_r \\ d_r & c_r & f_r \end{matrix}\right\}_{q=\xi_r^2}^{RW}
\right|
=
{\rm Vol}(T), 
\label{eq:volume}
\end{equation}
where ${\rm Vol}(T)$ is the hyperbolic volume of $T$.  
\end{theorem}
By specifying the sequences $a_r$, $b_r$, $\cdots$, $f_r$ more
precisely, we also have the following refined asymptotics.  
\begin{theorem}
Let $T$ be a hyperbolic tetrahedron with at least one of of whose vertices are  ideal or ultra-ideal. 
Let  $\theta_a$, $\theta_b$, $\theta_c$, $\theta_d$, $\theta_e$, $\theta_f$ be dihedral angles of $T$
and
$\eta_a = \frac{\pi-\theta_a}{2\pi}$, $\cdots$, $\eta_f = \frac{\pi-\theta_f}{2\pi}$.    
Assume that $\eta_a$, $\cdots$, $\eta_f$ are all rational numbers whose denominators and numerators are odd integers.  
%\begin{multline*}
%\left|\frac{2a_r+ 1}{r} -\frac{\pi - \theta_a}{2\, \pi}\right| = \min_a\left|\frac{2a+ 1}{r} -\frac{\pi - \theta_a}{2\, \pi}\right|,
%\\ \left|\frac{2b_r+ 1}{r} -\frac{\pi - \theta_b}{2\, \pi}\right| = \min_b\left|\frac{2b+ 1}{r} -\frac{\pi - \theta_b}{2\, \pi}\right|,
%\cdots,\ \  
%\left|\frac{2f_r+ 1}{r} -\frac{\pi - \theta_a}{2\, \pi}\right| = \min_f\left|\frac{2f+ 1}{r} -\frac{\pi - \theta_a}{2\, \pi}\right|.
%\end{multline*}
Let $r$ be an odd integer which is a common multiple of  the denominators of $\eta_a$, $\cdots$, $\eta_f$.
Then  $a_r = \frac{r}2 \eta_a - \frac{1}{2}$, $\cdots$, $f_r =  \frac{r}2 \eta_f - \frac{1}{2}$ are  integers  
and
\begin{equation}
\left|\left\{\begin{matrix}
a_r & b_r& e_r
\\ 
d_r & c_r & f_r
\end{matrix}\right\}_{q=\xi_{r}^2}^{RW}\right|
\ \ 
\underset{\begin{matrix}\scriptstyle r\to\infty
%\\[-3pt]
%\scriptstyle \text{$r$ is odd}
\end{matrix}}\sim
\ \ 
\frac{\sqrt{2} \, \pi}{{r}{}^{3/2} \, \sqrt[4]{-\det G}} \, e^{\frac{{r}}{2\, \pi}\, {\rm Vol}(T)}, 
\label{eq:th}
\end{equation}
where $G$ is the Gram matrix of $T$ given by
\[
G = 
\begin{pmatrix}
1 & -\cos \theta_a & -\cos \theta_b & -\cos\theta_f \\
-\cos\theta_a & 1 &  -\cos \theta_e & -\cos\theta_c \\
-\cos\theta_b & -\cos\theta_e & 1 & -\cos\theta_d \\
-\cos\theta_f & -\cos\theta_c & -\cos\theta_d & 1
\end{pmatrix}.  
\]
%Here we are considering  the asymptotics with respect to large $r$ where $r$ is odd and all of  $a_r$, $\cdots$, $f_r$ are integers.  
%}
\end{theorem}
\begin{remark}
For the classical $6j$ symbol, such asymptotics were conjectured by Wigner, Ponzano and Regge, and proved by J. Roberts \cite{R1}. 
For non-squared $q$  case,  the asymptotics of the quantum $6j$ symbol corresponding to  a spherical tetrahedron was given by Taylor and Woodward  in \cite{TW}.  
Here we use squared $q=\xi_r^2$, then the quantum $6j$ symbol corresponding to a hyperbolic tetrahedron grows exponentially and the leading term is given by the volume of the tetrahedron.   
This point is quite different from the asymptotics given in the previous works, but 
our asymptotics also relate to the determinant of the Gram matrix as in the previous works.  
\end{remark}
\par
The paper is organized as follows.  
In Section 2, we explain the quantum $6j$ symbol, and  
prove Theorem 1 in Section 3 that the leading term of the asymptotics of the quantum $6j$ symbol is given by the hyperbolic volume of the corresponding tetrahedron.  
%In Section 3,  a strategy to prove the volume conjecture is explained, which we use in the next section.   
We prove Theorem 2 in Section 4 by showing that  the second leading therm of the asymptotics of the quantum $6j$ symbols is given by the determinant of the corresponding Gram matrix.      
Some part of the proof of Theorem 3 is given in Appendix.  
In Section 4, a conjecture about the symmetry of the asymptotics of the quantum $6j$ symbols is proposed.  
\section{Quantum $6j$ symbols}
The quantum $6j$ symbol was introduced by Kirillov and Reshetikhin \cite{KR} 
to describe the structure of the tensor representations of the quantum group 
${\mathcal U}_q(sl_2)$.  
Let $q$ be an indeterminate.  
Finite-dimensional irreducible representations of ${\mathcal U}_q(sl_2)$ are
parametrized by non-negative half integers, which is called the spin.  
For a non-negative half integer  $a$, let $V_a$ denote the corresponding finite-dimensional irreducible representation.  
It is known that $V_c$ is contained in $V_a \otimes V_b$ if $a$, $b$, $c$ satisfied the Clebsch-Gordan condition, which means that $|a-b| \leq c \leq a+b$ and $a+b+c$ is an integer.  
Let $i_{ab}^c$ denote the inclusion map from $V_c$ to $V_a \otimes V_b$.  
For six non-negative half integers $a$, $b$, $c$, $d$, $e$, $f$, the quantum $6j$ symbol 
is defined by the following relation between the following two inclusions of an irreducible representation to a tensor of three irreducible representations $V_c \to V_d\otimes V_b \otimes V_a$, one through $V_d \otimes V_e$ and another one through $V_f \otimes V_a$.  
\begin{multline*}
\left(V_c \overset{i_{de}^c}\longrightarrow V_{d} \otimes V_{e} 
\overset{id \otimes i_{ba}^e}\longrightarrow V_{d} \otimes V_b \otimes V_a\right)
=
\\
\sum_f
\left\{\begin{matrix}
a & b & e \\ d & c & f\end{matrix}\right\}_q\,
\left(V_c \overset{i_{fa}^c}\longrightarrow V_{f} \otimes V_{a} 
\overset{i_{db}^f\otimes id}\longrightarrow V_{d} \otimes V_b \otimes V_a\right).
\end{multline*}
Here we assume that the triplets $(a, b, e)$, $(a, c, f)$, $(b,d, f)$ and $(c, d, e)$ satisfy the Clebsch-Gordan condition.  
The Racah-Wigner version of the quantum $6j$ symbol is given by
$$
\left\{\begin{matrix}
a & b & e \\ d & c & f\end{matrix}\right\}_q^{RW}
=
\frac{1}{\sqrt{-1}^{2(c+d+2e-a-b)}\,\sqrt{[2\, e+1]\, [2\, f+1]}}
\left\{\begin{matrix}
a & b & e \\ d & c & f\end{matrix}\right\}_q,  
$$
where $[n] = \frac{q^{n/2}-q^{-{n/2}}}{q^{1/2}-q^{-1/2}}$.  
\textcolor{black}{
Let $[n]! = \prod_{j=1}^n [j]$.  
}
Then this version has the following symmetric formula
\begin{multline}
\left\{\begin{matrix}
a & b & e \\ d & c & f\end{matrix}\right\}_q^{RW}
=
\\
\Delta_r(a, b, e)\, \Delta_r(a, c, f) \, \Delta_r(b, d, f) \, \Delta_r(c, d, e) %\, \times
%\\
\sum_{z=m}^{M}
\frac{(-1)^z\, [z+1]!}
{\prod_{j=1}^4[z-S_j]! \,\prod_{k=1}^3 [T_k-z]! },
\label{eq:sixj}
\end{multline}
where
\begin{equation}
\begin{aligned}
S_1& = a+b+e, \quad S_2 = a+c+f, \quad S_3 = b+d+f, \quad S_4 = c+d+e,
\\
T_1 &= a+b+c+d, \quad T_2 = a + d+e+f, \quad T_3 = b+c+e+f,
\\
m &=\max(S_1,S_2, S_3, S_4), \quad M =\min(T_1, T_2, T_3),
\\
\Delta_r&(u, v, w) = \left(
\frac{[u+v-w]! \, [v+w-u]! \, [w+u-v]!}{[u+v+w+1]!}\right)^{1/2}.  
\end{aligned}
\label{eq:ST}
\end{equation}
\par
For quantum invariants of knots, its relation to the hyperbolic volume of the knot complement was found by R. Kashaev \cite{K}, and reformulated  using the colored Jones polynomial in \cite{MM}.  
This relation is called the volume conjecture and proved for a few knots and links.    
The volume conjecture suggests that there must be a relation between the quantum $6j$ symbol and the volume of a hyperbolic tetrahedron.  
Using this idea, a volume formula for hyperbolic and spherical tetrahedra was obtained by the second author and M. Yano in \cite{MY}, and  the
asymptotics of the quantum $6j$ symbol was studied by F. Costantino \cite{C}.
\par
Let
\[
\left\{\begin{matrix}
a & b & e \\ d & c & f\end{matrix}\right\}_q^{RW}
\]
be the Racah-Wigner version of the quantum $6j$ symbol which is defined for six non-negative half integers $a$, $b$, $\cdots$, $f$ with the quantum parameter $q$.  
A triplet $(a, b, c)$ of half integers is called {\it admissible} if it satisfies the Clebsch-Gordan condition, i.e.  $|a-b| \leq c \leq a+b$ and $a+b+c \in {\mathbf Z}$.  
\begin{theorem}[F. Costantino \cite{C}]
Let $T$ be a hyperbolic truncated tetrahedron with dihedral angles $\theta_a$, $\theta_b$, $\theta_c$, $\theta_d$, $\theta_e$, $\theta_f$ at edges $a$, $\cdots$, $f$ in Figure \ref{fig:tet}. 
Let  $a_r$, $b_r$, $\cdots$, $f_r$ be sequences of non-negative half integers satisfying 
$$
2\, \pi \, \left|\lim_{r\to\infty} \frac{a_r}{r} - \frac{1}{2}\right| = \theta_a, \ 
2\, \pi \, \left|\lim_{r\to\infty} \frac{b_r}{r} - \frac{1}{2}\right| = \theta_b, \   
\cdots, \ 
2\, \pi \, \left|\lim_{r\to\infty} \frac{f_r}{r} - \frac{1}{2}\right| = \theta_f,
$$
and the triplets $(a_r, b_r, e_r)$, $(a_r, c_r, f_r)$, $(b_r, d_r, f_r)$ and  $(c_r, d_r, e_r)$ are all admissible% and of $r$-hyperbolic type for $r \in {\mathbb N}$
.  
Then
\[
\lim_{r\to\infty}\frac{2\, \pi}{r} \, \log\left|{\rm ev}_r\left\{\begin{matrix}
a_r & b_r & e_r \\ d_r & c_r & f_r
\end{matrix}\right\}_{\xi_r}^{RW}\right|
=
2 \, {\rm Vol}(T),  
\]
where $\xi_r = \exp(2 \pi \sqrt{-1}/r)$ and ${\rm Vol}(T)$ is the hyperbolic volume of \,$T$.  
\end{theorem}
Here the operator $\mathrm{ev}_r$ removes the pole of the quantum $6j$ symbol at $q = \xi_r$.  
In this paper, we put $q = \xi_r^2$ instead of $\xi_r$,  
then the quantum $6j$ symbol doesn't have a pole at $q=\xi_r^2$ and its actual asymptotics is given by the volume and the determinant of the Gram matrix of the corresponding tetrahedron.  
\section{Relation to the hyperbolic volume}
In this section, we prove Theorem 1 giving a relation between the Racah-Wigner version of the quantum 6j symbols and the hyperbolic volume of the corresponding hyperbolic tetrahedra.  
%
%$$
%\lim_{r\to\infty}2\pi\, \tfrac{2a_r+1}{r}  = \pi - \theta_a, \quad
%\lim_{r\to\infty} 2\pi \,\tfrac{2b_r+1}{r}  = \pi - \theta_b, \quad\cdots,\quad 
%\lim_{r\to\infty} 2\pi\, \tfrac{2f_r+1}{r}  = \pi - \theta_f
%$$  
%such that the triplets $(a_r, b_r, e_r)$, $(a_r, d_r, f_r)$, $(b_r, d_r, f_r)$,  $(c_r, d_r, e_r)$ are all $r$-admissible for odd $r \geq 3$.  
%Let $a_r=\frac{r}{2} \, \eta_a - \frac{1}{2}$, $\cdots$, $f_r = \frac{r}{2} \, \eta_f - \frac{1}{2}$.  
%Then
%\begin{equation}
%\lim_{r\to\infty}\frac{2\, \pi}{r} 
%\log\left|
%\left\{\begin{matrix}
%a_r & b_r & e_r \\ d_r & c_r & f_r \end{matrix}\right\}_{q=\xi_r^2}^{RW}
%\right|
%=
%{\rm Vol}(T), 
%\label{eq:volume}
%\end{equation}
%where ${\rm Vol}(T)$ is the hyperbolic volume of $T$.  
%Here we are considering the asymptotics with respect to large $r$ where $r$ is odd and all of $a_r$, $\cdots$, $f_r$ are integers.   
%\end{prop}
%
\par\noindent
{\it Proof of Theorem 1}. \ 
The idea of proof is similar to the proof of the volume conjecture for figure-eight knot, for example, see \cite{M}.
Such an idea is also used in \cite{C} for the quantum $6j$ symbols with $q=\xi_r$, where he used $\mathrm{ev}_r$ map.  
 Here we follow the argument in \cite{C} without $\mathrm{ev}_r$ map.  
At the end of proof, the limiting value is compared with the volume formula in \cite{U}  based on \cite{MY}.  
\par
From the assumption, 
$m = \max(S_{r,1}, \cdots, S_{r,4}) \geq \frac{r-3}{2}$, $M=\min(T_{r,1}, T_{r,2}, T_{r,3}) < \frac{3r}{4}$  and  $\frac{r-1}{2} \leq m+1 \leq z+1 \leq r-1$.  
For such $z$, $(-1)^z \, [z+1]!$ has the same sign since $[k+1] < 0$ for $m+1 \leq k \leq M$.   
Therefore the terms in the sum of \eqref{eq:sixj} are all real numbers with the same sign, and
the limit in \eqref{eq:volume} is determined by the term  having the largest absolute value.  
%}
Indeed, let
\begin{equation}
\alpha_r(z) =
\left.\frac{(-1)^z\, [z+1]!}
{\prod_{j=1}^4 [z-S_{r,j}]!  \,\prod_{k=1}^3 [T_{r,k}-z]!}\right|_{q=\xi_r^2},
\label{eq:alpha}
\end{equation}
where $S_{r,j}$, $T_{r,j}$ are defined by using $a_r$, $b_r$, $\cdots$, $f_r$
as in \eqref{eq:ST},
\begin{equation}
m_r = \max_{j=1}^4 S_{r,j}, \qquad M_r = \min_{k=1}^3 T_{r,k}
\label{eq;mM}
\end{equation}
and let $|\alpha_r(z_0)|$ be the maximum of $\{|\alpha_r(z)| \mid m_r \leq z \leq M_r\}$.  
Then
\[
\lim_{r\to\infty} \frac{1}{r} \log \left| \alpha_r(z_0)\right|
<
\lim_{r\to\infty} \frac{1}{r} \log \left| \sum_{z=m}^M \alpha_r(z)\right|
<
\lim_{r\to\infty} \frac{1}{r} \log \left|r\, \alpha_r(z_0)\right|.  
\]
Since $\displaystyle\lim_{r\to\infty} \frac{1}{r} \log r = 0$, we get 
$\displaystyle\lim_{r\to\infty} \frac{1}{r} \log \left| \sum_{z=m}^M \alpha_r(z)\right|
=
\displaystyle\lim_{r\to\infty} \frac{1}{r} \log \left|\alpha_r(z_0)\right|$.  
\par
Now compare 
\[
\lim_{r\to\infty} \frac{2 \, \pi}{r}
\log \left|
\Delta_r(a_r, b_r, e_r) \, \Delta_r(a_r, c_r, f_r) \, \Delta_r(b_r, d_r, f_r) \, \Delta_r(c_r, d_r, e_r) \, \alpha_r(z_0)\right|  
\]
with the hyperbolic volume of $T$.  
This part is the same as  the proof of Theorem 1 in \cite{C}, which concerns  the quantum $6j$ symbol with $q=\xi_r$.  
Here we assume $r$ to be odd and $q$ to be $\xi_r^2$ instead of $\xi_r$, so the proof in \cite{C} works well without applying the operator $\mathrm{ev}_r$ to remove the factor $\xi_r^{r/2} - \xi_r^{-r/2}$ in $[S_{r,j}+1]!$ and $[z+1]!$.  
Here we may have factors $\xi_r^{(r\pm1)/2} - \xi_r^{-(r\pm1)/2}$ in  $[S_{r,j}+1]!$ and $[z+1]!$ but we never have the factor $\xi_r^{r/2} - \xi_r^{-r/2}$ since $q = \xi_r^2$ and $r$ is odd.  \qed
%\end{proof}
%
\medskip
\par\noindent
Theorem 1 can suggest the following a bit more general conjecture which is also inspired by the new volume conjecture proposed by the first author and T. Yang \cite{CY}. 
\begin{conjecture}
Let $T$ be a hyperbolic tetrahedron and
let  $\theta_a$, $\theta_b$, $\theta_c$, $\theta_d$, $\theta_e$, $\theta_f$ are dihedral angles of $T$
and
$a_r$, $b_r$, $\cdots$, $f_r$ be sequences of admissible half integers for odd $r\geq3$  satisfying \eqref{eq:sequence} as in Theorem 1.
%$$
%\lim_{r\to\infty} \frac{a_r}{r} = \frac{1}{4}-\frac{\theta_a}{4\pi}, \ 
%\lim_{r\to\infty} \frac{b_r}{r} = \frac{1}{4}-\frac{\theta_b}{4\pi}, \ 
%\cdots, \ 
%\lim_{r\to\infty} \frac{f_r}{r} = \frac{1}{4}-\frac{\theta_f}{4\pi} .
%$$
Then
\begin{equation}
\lim_{r\to\infty}\frac{2\, \pi}{r} 
\log\left|
\left\{\begin{matrix}
a_r & b_r & e_r \\ d_r & c_r & f_r \end{matrix}\right\}_{q=\xi_r^2}^{RW}
\right|
=
{\rm Vol}(T).
\label{eq:volumeconjecture}
\end{equation}
\end{conjecture}
\section{Relation to the Gram matrix}
In this section, we investigate a bit more about the asymptotics of $\alpha_r(z)$ and prove Theorem~2 by using the method to prove the volume conjecture for $5_2$ knot given in \cite{O}.   
Here we use notations
%\[
%a_r = \frac{r\, \eta_a-1}{2}, \quad
%b_r = \frac{r\, \eta_b-1}{2}, \quad \cdots, \quad
%f_r = \frac{r\, \eta_f-1}{2},
%\]
%and 
$S_{r,j}$, $T_{r,k}$,
$m_r $, $M_r$  introduced in \eqref{eq:alpha}.  
%\medskip
\par
%{\bf Step 1.}
First of all, we express $a_r(z)$ in terms of the dilogarithm function.  
Let  
\[
(q)_n = \prod_{k=1}^n (1-q^k).
\]  
Then 
\[
[n]! = 
\textcolor{black}{
q^{\frac{-n(n+1)}{4}}\,\frac{(q)_n}{2^n \,(-i)^n \,  \sin^n \frac{2\pi}{r}}
}
\]
and $\alpha_r(z)$, $\Delta_r(a_r, b_r, e_r)$, $\cdots$, $\Delta_r(c_r, d_r, e_r)$ are reformulated as
\begin{equation}
\alpha_r(z)= 
\frac{ \xi_r^{2(d_1+d_2)}\, (-1)^{z+1} \, (\xi_r^2)_{z+1}}{2\, i \, \sin\frac{2\pi}{r}\,\prod_{j=1}^4 (\xi_r^2)_{z-S_j} \,  \, \prod_{k=1}^3(\xi_r^2)_{T_k-z} },
\label{eq:6jq}
\end{equation}
where 
\[
\begin{aligned}
d_1 &= 
%- ( i^2 + j^2 + k^2 + l^2 + m^2 +n^2+ i j+
%   i k+ i l+ i m+i n+ 
%   j k+ j l+ j m+ j n+ 
%   \\
%   k l+ k m+ k n+ l m+
%   l n+ m n)+
-\left(2a_r+2 b_r+2 c_r+2 d_r+2 e_r+2
   f_r + \frac{1}{2}\right)z+\frac{3 z^2}{2} ,
\\
d_2 &= 
a_r^2 + a_r b_r + b_r^2 + a_r c_r + b_r c_r + c_r^2 + a_r d_r + b_r d_r + c_r d_r + d_r^2 + a_r e_r + 
\\&
\qquad
 b_r e_r + c_r e_r + d_r e_r + e_r^2 + a_r f_r + b_r f_r + c_r f_r + d_r f_r + e_r f_r + f_r^2 - \frac{1}{2},
 \end{aligned}
\]
and
\begin{align*}
\Delta_r(a_r, b_r, e_r)
&=
\left(
\xi_r^{2d_3(a_r, b_r, e_r)} \,
\frac{
\textcolor{black}{-}
\,
2 \, i \, \sin\frac{2 \pi}{r} \, (\xi_r^2)_{a_r+b_r-e_r} \, (\xi_r^2)_{b_r+e_r-a_r} \, (\xi_r^2)_{e_r+a_r-b_r}}{(\xi_r^2)_{a_r+b_r+e_r+1}}
\right)^{1/2},
\\
&\cdots,
\\
\Delta_r(c_r, d_r, e_r)
&=
\left(
\xi_r^{2d_3(c_r, d_r, e_r)} \,
\frac{
\textcolor{black}{-}
\,
2 \, i \, \sin\frac{2 \pi}{r} \, (\xi_r^2)_{c_r+d_r-e_r} \, (\xi_r^2)_{d_r+e_r-c_r} \, (\xi_r^2)_{e_r+c_r-d_r}}{(\xi_r^2)_{c_r+d_r+e_r+1}}
\right)^{1/2}
\end{align*}
where 
\[
d_3(a_r, b_r, e_r) = -\frac{1}{2} \, (a_r^2+b_r^2+e_r^2-2a_rb_r-2a_re_r-2b_re_r-a_r-b_r-e_r-1).
\]  
\par
%The signature of $\alpha_z^{(r)}$ is $(-1)^{(r+3)/2}$.  
Let $\varphi_r$ is the function introduced in \cite{F} for $0 \leq t \leq 1$ by the following integral.  
\[
\varphi_r(t)
=
\int_{-\infty}^{\infty}
\frac{e^{(2t-1)x}}{4\,x \, \sinh x \sinh(2x/r)} \, dx, 
\]
where the singularity at $x=0$ is bypassed by a small half circle above the real line.  
\begin{lemma}
The function $\varphi_r$ satisfies
\begin{align}
(1-e^{2\pi i a}) 
&=
\exp\left(\varphi_r\left(a-\frac{1}{r}\right) 
- 
\varphi_r\left(a+\frac{1}{r}\right)\right)  &\left(\frac{1}{r} \leq a \leq \frac{r-1}{r}\right),
\label{eq:res}
\\
(\xi_r^2)_n &=
\exp\left(\varphi_r\left(\frac{1}{r}\right) 
- 
\varphi_r\left(\frac{2n+1}{r}\right)\right)  &\left(0 \leq n \leq \frac{r-1}{2}\right), 
\label{eq:qn1}
\\
(\xi_r^2)_n &= 
\exp\left(\varphi_r\left(\frac{1}{r}\right) 
- \varphi_r\left(\frac{2n+1}{r}-1\right)+\log 2 \right) 
&\left(\frac{r-1}{2}<n < \frac{2r-1}{2}\right),
\label{eq:qn2}
\\
{\rm Im}\,\varphi_r(t) &= -\pi \, i \, \frac{6 \,r^2 \,t^2-6\, r^2\,  t+r^2-2}{24\, r}.  &{ }
  \label{eq:Im}
\end{align}
\end{lemma}
\begin{proof}
We have 
$$
\begin{aligned}
\varphi_r(a-\frac{1}{r})-\varphi_r(a+\frac{1}{r})  
&=
\int_{-\infty}^\infty \frac{e^{(2a-1-\frac{2}{r})x}-e^{(2a-1+\frac{2}{r})x}}{4\,x \, \sinh x \sinh(2x/r)} \, dx
\\
&=
-\int_{-\infty}^\infty \frac{e^{(2a-1)x}}{4\,x \, \sinh x } \, dx.
\end{aligned}
$$
The function $-\frac{e^{(2a-1)x}}{2x \, \sinh x}$ has poles at $x=k\pi i$ $(k \in \mathbb N)$ in the upper half plane, and the corresponding residues are $-\frac{e^{2a\pi k i}}{2 \pi k i}$.  The sum of these residues is equal to $\frac{\log(1-e^{2a\pi  i})}{2 \pi i}$ and this implies \eqref{eq:res}. 
\par
The formula \eqref{eq:qn1} is a product of \eqref{eq:res} for $a=2/r$, $4/r$, $\cdots$, $2n/r$.  
\par
To prove \eqref{eq:qn2}, we decompose 
$
(\xi_r^2)_n$
into a product of 
$
(\xi_r^2)_{\frac{r-1}{2}}
$
and 
$
\prod_{k=1}^{n-\frac{r-1}{2}} (1 - \xi_r^{2k-1})
$.  
Then $(\xi_r^2)_{\frac{r-1}{2}} = \exp(\varphi_r(1/r)-\varphi_r(1))$ and 
$\prod_{k=1}^{n-\frac{r-1}{2}} (1 - \xi_r^{2k-1}) = \exp(\varphi_r(0)-\varphi_r(\frac{2n+1}{r}-1))$.  
Moreover, $\varphi_r(0) - \varphi_r(1)$ is given by $2 \pi  i$ times the sum of residues of $-\frac{1}{2\,x \,\sinh(2x/r)}$, which is   $\sum_{k=1}^\infty \frac{(-1)^{k+1}}{ k} = \log 2$.  
Hence we get \eqref{eq:qn2}.  
\par
The imaginary part of $\varphi_r(t)$ comes from the pole at $x=0$.  
For the other part, $\varphi_r(t)$ is defined as an integral of a real function.  
We define the integral path to avoid $x=0$ by taking a path in the upper half plane, so by taking a small upper half circle around $x=0$, the contribution is given by $- \pi \, i$ times the residue of
$\frac{e^{(2t-1)x}}{4\,x \, \sinh x \sinh(2x/r)}$ at
$x=0$, which is the right hand side of \eqref{eq:Im}.  
\end{proof}
Now, we reformulate $\alpha_r(z)$ and $\Delta_r(a_r,b_r,e_r)$ as continuous functions by using \eqref{eq:qn1} and \eqref{eq:qn2}.  
%Let
%$m_r = \max_{1 \leq j \leq 4} S_{r,j}$, $M_r = \min_{1 \leq k \leq 3}T_{r,k}$, 
\begin{multline}
\widetilde\alpha_r(z) = 
%\frac{r}{4\, \pi} \, \times
%\\
i^r \, \xi_r^{2d_2}
 \, \exp\left(
%\textcolor{red}{ - }\,
-\pi\, i \, z - 6\, \varphi_r(\tfrac{1}{r}) - \varphi_r(\tfrac{2z+3}{r}-1)+ \log 2
+ \sum_{j=1}^4\varphi_r(\tfrac{2z-2S_{r,j}+1}{r})
\right.
\\
\left.
+ \sum_{j=1}^3\varphi_r(\tfrac{2T_{r,j}-2z+1}{r})
+ \frac{2 \pi i}{r}\left(3z^2 - z - 4\,(a_r + b_r + \cdots + f_r)\,z \right)
\right), 
\label{eq:contalpha}
\end{multline}
and
\begin{multline}
\widetilde\Delta_r(a_r, b_r, e_r) 
=
\exp\left(
\frac{1}{2} \, \left( 
2 \, \varphi_r(\tfrac{1}{r})+
\varphi_r(\tfrac{2a_r+2b_r+2e_r+3}{r}-1) - 
\right.\right.
\\
\left.\left.
\log 2-\varphi_r(\tfrac{2a_r +2b_r-2e_r+1}{r}) -
\varphi_r(\tfrac{2b_r+2e_r-2a_r+1}{r}) - 
\varphi_r(\tfrac{2e_r+2a_r-2b_r+1}{r})
\vphantom{\frac{4 \pi i}{r}}
\right)\right). 
\label{eq:contdelta}
\end{multline}
Then 
\textcolor{black}{
$\tilde\alpha_r(z)$ is a positive real function for $m_r< z < M_r$, 
}
\[
\alpha_r(z) = \frac{\textcolor{black}{
(-1)^{\frac{r+1}{2}} 
}
}{2\,
 \sin\frac{2\pi}{r}}\,\widetilde \alpha_r(z)
 \] 
 for an integer $z$ $(m_r \leq z \leq M_r)$ and   
 \[
 \Delta_r(a_r, b_r, e_r) 
 = 
\sqrt{-2\, i\,  \sin\frac{2\pi}{r}}\,\xi_r^{d_3(a_r,b_r,e_r)}\,
 %\textcolor{black}{i^{-\frac{r-3}{2}+(a_r+b_r+e_r)}\,}
 \widetilde\Delta_r(a_r, b_r, e_r).
 \] 
% for an $r$-admissible triplet $(a_r, b_r, e_r)$.  
%\medskip
\par
%{\bf Step 2.} 
Next, we apply the Poisson summation formula to $\sum_{z=m_r}^{M_r} \widetilde\alpha_r(z)$.  
Let $\varepsilon$ be a small positive real number and $\psi(x)$ be a  smooth function satisfying $0 \leq \psi(x) \leq 1$, $\psi(0)=0$ for $x\in \mathbb R \setminus \left((2m_r+3)/r+ \varepsilon, (2M_r+3)/r-\varepsilon\right)$, and $\psi(x)=1$ for $x \in [(2m_r+3)/r +2\varepsilon, (2M_r+3)-2\varepsilon]$. 
Note that   $(2m_r+3)/r$ and $(2M_r+3)/r$ do not depend on $r$.  
Let 
\begin{equation}
h_r(z) = \psi\left(\frac{2z+3}{r}\right) \, \widetilde \alpha_r(z).
\label{eq:psi}
\end{equation} 
\par
By the Poisson summation formula, we have
\[
\sum_{m\in {\mathbf Z}}h_r(m)
= 
\sum_{m\in {\mathbf Z}} \widehat h_r(m)
\]
where 
$\widehat h_r(m) = \int_{-\infty}^\infty e^{-2 \pi i m t}\, h_r(t)\, dt$.  
The function $h_r(t)$ is non-negative or non-positive depending on $r$, and satisfies the following.  
\begin{lemma}
There are positive constants $c$ and $\delta'$ satisfying
\begin{equation}
\frac{\left|\sum_{m\neq1}\widehat h_r(m)\right|}{\left|\widehat h_r(0)\right|} < 
c \, e^{-r\, \delta}.  
\label{eq:limit}
\end{equation}
\label{lem:remain2}
\end{lemma}
We also have the following.  
\begin{lemma}
Let $A_r = \sum_{z=m_r}^{M_r} \widehat\alpha_r(z)$ and $B_r = \sum_{z=m_r}^{M_r} h_r(z)$, then there is positive constants $c'$ and $\delta'$ not depending on $r$ satisfying
\begin{equation}
\frac{|A_r - B_r|}{|\widehat h_r(0)|} < c' \, e^{-\delta'}.  
\label{eq:remain}
\end{equation}
\label{lem:remain}
\end{lemma}
The proof of the above two lemmas are given in Appendix.  
\par
The above two lemmas shows that the main asymptotics of the quantum $6j$ symbol is determined by the asymptotics of $\widehat h_r(0)$ with respect to $r$.  
Let $\zeta = \frac{2z+3}{r}$, 
%$\eta_{a_r} = \frac{2 \, a_r+1}{r}$, 
%$\eta_{b_r} = \frac{2 \, b_r+1}{r}$, $\cdots$, $\eta_{ f_r} = \frac{2 \, f_r+1}{r}$, 
$\sigma_{j} = \frac{2\, S_{r,j} + 3}{r}$ $(1 \leq j \leq 4)$,  
and
$\tau_{k} = \frac{2\, T_{r,k}+4}{r}$ $(1 \leq k \leq 3)$. 
Note that $\sigma_j$, $\tau_k$ do not depend on $r$ since they are determined by $\eta_a$, $\cdots$, $\eta_f$.  
Let
\begin{multline}
g_r(\zeta, \boldeta) = 
%\frac{r}{4\, \pi} \, 
\\
\psi(\zeta) \, 
%\\
\exp\left(
%\textcolor{red}{-}\,
\frac{-\pi i r  \zeta}{2} + \frac{3 \pi i}{2}
- 6\, \varphi_r(\tfrac{1}{r}) - \varphi_r(\zeta-1)+ \log 2
+ \sum_{j=1}^4\varphi_r(\zeta-\sigma_{r,j}+\tfrac{1}{r})
\right.
\\
\left.
+ \sum_{j=1}^3\varphi_r(\tau_{r,j}-\zeta)
+ 2 \pi i r\left(\frac{3}{4}\zeta^2  - (\eta_{a_r} +  \cdots + \eta_{f_r})\,\zeta +\frac{\zeta}{r} \right)
\right) ,
\label{eq:g}
\end{multline}
where $\boldeta=(\eta_a, \eta_b, \cdots, \eta_f)$.  
Then the function
$h_r(z)$ is equal to 
$\left|g_r(\frac{2z+3}{r}, \boldeta)\right|$, 
the argument of $g(\zeta, \boldeta)$ does not depend on $\zeta$
and
\[
\left|\widehat h_r(0)\right| = 
\left|\int_{-\infty}^\infty h_r(x)\, dx\right| = 
\frac{r}{2} \left|\int_{-\infty}^{\infty}g_r(\zeta, \boldeta) \, d\zeta  \right|.
\]
\par
%Since $(a_{r_i}, \cdots, f_{r_i})$ is an infinite adequate sequence, we have
%\[
%\eta_{a_{r_i}} = \frac{1}{2\pi}\, (\pi - \theta_a) + O(\frac{1}{{r_i}^2}), 
%\ \cdots, \ 
%\eta_{f_{r_i}} = \frac{1}{2\pi}\, (\pi - \theta_f)+ O(\frac{1}{{r_i}^2}).
%\]
%Now we put
%\[
%\eta_{a} = \frac{1}{2\pi}\, (\pi - \theta_a), 
%\ \cdots, \ 
%\eta_{f} = \frac{1}{2\pi}\, (\pi - \theta_f),
%\]
%%
%\begin{multline}
%\bar g_r(\zeta) = 
%\psi_r(\tfrac{r \zeta-3}{2}) \, \times
%\\
%\exp\left(
%\frac{\pi i r  \zeta}{2} - 6\, \varphi_r(\tfrac{1}{r}) - \varphi_r(\zeta-1)+ \log 2
%+ \sum_{j=1}^4\varphi_r(\zeta-\sigma_{j}+\tfrac{1}{r})
%\right.
%\\
%\left.
%+ \sum_{j=1}^3\varphi_r(\tau_{j}-\zeta)
%+ 2 \pi i r\left(\frac{3}{4}\zeta^2  - (\eta_{a} +  \cdots + \eta_{f})\,\zeta +\frac{\zeta}{r}\right)
%\right).  
%\label{eq:barg}
%\end{multline}
%%
%Then
%\[
%g_{r_i}(\zeta) = \bar g_{r_i}(\zeta)\times \exp(O(\frac{1}{r_i})).  
%\]
%
The parameters $\eta_a$, $\cdots$, $\eta_f$ do not depend on $r$, but the function $\varphi_r(x)$ depends on $r$.  
We investigate the asymptotics of 
\begin{equation}
\left|r \sin\tfrac{2 \pi}{r} \,\widetilde \Delta(\tfrac{r \eta_a -1}{2}, 
\tfrac{r \eta_b -1}{2}, \tfrac{r \eta_e -1}{2}) \cdots \widetilde \Delta(\tfrac{r \eta_c -1}{2}, \tfrac{r \eta_d -1}{2}, \tfrac{r \eta_e -1}{2}) \, \int_{-\infty}^{\infty} g_r(\zeta, \boldeta)\, d\zeta\right|.
\label{eq:recont}
\end{equation}
We expand $g_r(\zeta, \boldeta)$ %in terms of $r$
by using the following relations.  
%These properties are explained in Appendix A of \cite{O}.  
%
\begin{lemma}[{\cite[Appendix A]{O}}]
The asymptotics of $\varphi_r$ is given as follows.  
\begin{align}
\varphi_r\left(\tfrac{1}{r}\right) &= 
\frac{r}{4 \pi i} \left(
\frac{\pi^2}{6} + \frac{2 \pi i}{r} \log\frac{r}{2} - \frac{\pi^2}{r} + 
O(\tfrac{1}{r^2})\right),
\label{eq:r1}
\\
\varphi_r(t) &=
\frac{r}{4 \pi i}\left({\rm Li}_2(e^{2 \pi i t}) + O(\tfrac{1}{r^2})\right).  
\label{eq:r2}
\end{align}
\end{lemma}
We also have the following expansion for small $\alpha$ since $\frac{d}{dx}{\Li(x)}=-\frac{log(1-x)}{x}$.  
\begin{equation}
{\rm Li}_2(e^{2 \pi i(p+\alpha)}) =
{\rm Li}_2(e^{2 \pi i p}) - 2 \, \pi i \log(1-e^{2 \pi i p}) \, \alpha + O(\alpha^2). 
\label{eq:Li}
\end{equation}
By reformulating $g_r$ by using \eqref{eq:r1}, \eqref{eq:r2} and \eqref{eq:Li}, we get
\begin{equation}
\begin{aligned}
g_r&(\zeta, \boldeta) 
= 
\frac{16 \, \psi(\zeta)}{r^3}\, \exp\left(\!
\frac{-\pi i  r \zeta}{2}   + \frac{3 \pi i}{2} + \right.
\\
&
  \frac{r}{4 \pi i}\left(\!-\pi^2  + \frac{6\pi^2}{r}-{\rm Li}_2(e^{2 \pi i \zeta}) 
+ \!\sum_{j=1}^4{\rm Li}_2(e^{2 \pi i (\zeta-\sigma_j)}) + \!\sum_{j=1}^3{\rm Li}_2(e^{2\pi i(\tau_j-\zeta)} )
\right)+
\\
&\left.
2\pi i r\left(\frac{3}{4}  \zeta^2  -  (\eta_{a} + \cdots + \eta_{ f}){\zeta}+ \frac{\zeta}{r}\right) - 
 \frac{1}{2}  \sum_{j=1}^4\log( 1 - e^{2\pi i(\zeta-\sigma_j)}) + O(\tfrac{1}{r})\!
\right)
\\
=&
\tfrac{16 \, \psi(\zeta)}{r^3}\, 
\tfrac{-e^{2 \pi i \zeta}}{\sqrt{\prod_{j=1}^4 (1 - \tfrac{e^{2 \pi i \zeta}}{e^{2 \pi i \sigma_j}})}}\,%\times \qquad\qquad\qquad\qquad\qquad
%\\
%&\qquad
\exp\left(\!
\tfrac{r}{4 \pi i}\!\left(-\pi^2  -
2\pi^2 \zeta-{\rm Li}_2(e^{2 \pi i \zeta}) +
\right.
\right.
\\
&\left.
\left. 
\sum_{j=1}^4{\rm Li}_2(\tfrac{e^{2 \pi i \zeta}}{e^{2 \pi i \sigma_j}})  
+\sum_{j=1}^3{\rm Li}_2(\tfrac{e^{2 \pi i \tau_j}}{e^{2 \pi i \zeta}} )
-6   \pi^2 \zeta^2  + 8\pi^2\,(\eta_{a} + \cdots + \eta_{f})\,{\zeta}\right)
+ O(\tfrac{1}{r})
\right). 
\end{aligned}
\label{eq:gr}
\end{equation}
%
%
%We also reformulate $\widetilde\Delta_r$. 
Moreover, let
\begin{multline}
\delta(\eta_1, \eta_2, \eta_3) = 
\\
\frac{1}{8\pi i}
\left(\vphantom{\tfrac{e^{2 \pi i (\eta_1+\eta_2)}}{e^{2 \pi i \eta_3}}}
 {\rm Li}_2(e^{2 \pi i (\eta_1+\eta_2+\eta_3)})
- 
{\rm Li}_2(\tfrac{e^{2 \pi i (\eta_1+\eta_2)}}{e^{2 \pi i \eta_3}})- 
{\rm Li}_2(\tfrac{e^{2 \pi i (\eta_2+\eta_3)}}{e^{2 \pi i \eta_1}})
-{\rm Li}_2(\tfrac{e^{2 \pi i (\eta_3+\eta_1)}}{e^{2 \pi i \eta_2}})
\right)
\end{multline}
for $0 \leq \eta_i \leq 1$ and $1 \leq \eta_1+\eta_2+\eta_3 \leq 2$.  
Then
\begin{equation}
\widetilde\Delta_r(\tfrac{r\eta_1-1}{2}, \tfrac{r\eta_2-1}{2}, \tfrac{r\eta_3-1}{2}) 
= 
\frac{\sqrt{r}}{2}\, 
\,\exp\left(r\, \delta(\eta_1, \eta_2, \eta_3)
+ O(\tfrac{1}{r})\right).
\label{eq:adelta}
\end{equation}
\medskip\par
%{\bf Step 3.}  
Now we apply the stationary phase method to the function $g_r(\zeta, \boldeta)$.  
The arguments of 
\[
\tfrac{16 \, \psi(\zeta)}{r^3}
\,\tfrac{-e^{2 \pi i \zeta}}{\sqrt{\prod_{j=1}^4 (1 - \tfrac{e^{2 \pi i \zeta}}{e^{2 \pi i \sigma_j}})}}
\]
does not depend on $\zeta$ since
$\arg(1-z) = -(\pi-\arg z)/2$ for $z\in \mathbb C$ satisfying $|z|=1$.  
Similarly,  the argument of 
\begin{multline*}
\exp\left(
\frac{r}{4 \pi i}( -
2\pi^2 \zeta-{\rm Li}_2(e^{2 \pi i \zeta}) +
\sum_{j=1}^4{\rm Li}_2(\tfrac{e^{2 \pi i \zeta}}{e^{2 \pi i \sigma_j}})  
+
\right.
\\
\left.
\sum_{j=1}^3{\rm Li}_2(\tfrac{e^{2 \pi i \tau_j}}{e^{2 \pi i \zeta}} )
-6   \pi^2 \zeta^2  + 8\pi^2\,(\eta_{a} + \cdots + \eta_{f})\,{\zeta})
\right)
\end{multline*}
does not depend on $\zeta$ since 
$
\operatorname{Re} \,\mathrm{Li}_2(e^{2\pi x}) = \pi^2 ((x-\lfloor x\rfloor)^2-\frac{1}{2})^2 - \frac{\pi^2}{12}$
for $x \in \mathbb R
$
 where $\lfloor x \rfloor$ is the largest integer not exceeding $x$.  
 Let
\begin{multline*}
F(\zeta, \boldeta)
=
\frac{1}{4\, \pi \, i} \,
\left(%\frac{\pi^2}{6}  - \frac{\pi^2}{r}
\textcolor{black}{-}\,
2\pi^2 \zeta-{\rm Li}_2(e^{2 \pi i \zeta}) 
+\sum_{j=1}^4{\rm Li}_2(\tfrac{e^{2 \pi i \zeta}}{e^{2 \pi i \sigma_j}}) + 
\right.
\\
\left.
\sum_{j=1}^3{\rm Li}_2(\tfrac{e^{2 \pi i \tau_j}}{e^{2 \pi i \zeta}} )
-6   \pi^2 \zeta^2  + 8\pi^2\,(\eta_{a} + \cdots + \eta_{f})\,{\zeta}\vphantom{\sum_{j=1}^4}\right).
%\label{eq:gr}
\end{multline*}
\par
Let $\zeta_0$ $(m_r \leq \zeta_0 \leq M_r)$ be a solution of
\[
\frac{\partial}{\partial\zeta}
\mathrm{Re}\, F(\zeta, \boldeta)
=
0.
\]
 Then the derivative of $\operatorname{Im} F(\zeta, \boldeta)$ with respect to $\zeta$ is $0$ and $\zeta_0$ satisfies
\begin{equation}
\frac{\partial}{\partial\zeta}
 F(\zeta, \boldeta)
=
0.
\label{eq:quad}
\end{equation}
The equation \eqref{eq:quad} is essentially a quadratic equation and we have
only one solution which corresponds to the maximum. 
If we had two solutions corresponding to  maxima , we also should have at least one 
solution corresponding to a minimum, but the total number of the solutions of \eqref{eq:quad}
is two.  
\par
Since $\zeta_0 \in {\mathbf R}$ and the integral path is  the real line,  
we can apply the stationary 
phase method to $\left| \int_{-\infty}^\infty g_r(\zeta, \boldeta) \, d\zeta\right|
=\int_{-\infty}^\infty |g_r(\zeta, \boldeta)| \, d\zeta
$ and 
we have
\begin{equation}
\left| \int_{-\infty}^\infty g_r(\zeta, \boldeta) \, d\zeta\right|
\underset{r\to \infty}{\sim}
%\\
\frac{16}{r^3}\left|
\tfrac{e^{2 \pi i \zeta_0}}{{\prod_{j=1}^4}\sqrt{1-\tfrac{e^{2 \pi i \zeta_0}}{e^{2 \pi i \sigma_j}}}} \, 
\sqrt{\frac{2\, \pi}{r \left|F_{\zeta\zeta}(\zeta_0, \boldeta)\right|}} \, 
\exp\left(r F(\zeta_0, \boldeta)
+ O(\tfrac{1}{r})
\right)\right|.
\label{eq:ag}
\end{equation}
\begin{lemma}
\[
\frac{\prod_{j=1}^4 (1-\frac{e^{2 \pi i \zeta_0}}{e^{2 \pi i \sigma_j}})}{e^{4 \pi i \zeta_0}} \, F_{ss}(\zeta_0, \boldeta) 
= 
\pm4\, \pi \,  \sqrt{-\det G}
\]
where $G$ is the gram matrix of $T$.  
\end{lemma}
\begin{proof}
The second derivative of $F(\zeta, \boldeta)$ at $\zeta_0$ is given by 
\[
F_{\zeta\zeta}(\zeta_0, \boldeta)
=
\pi\, i  \, \left(
-\frac{e^{2 \pi i \zeta_0}}{1-e^{2 \pi i \zeta_0}} + 
\sum_{j=1}^4\frac{e^{2 \pi i \zeta_0}}{e^{2 \pi i \sigma_j}-e^{2 \pi i \zeta_0}}
-
\sum_{j=1}^3\frac{e^{2 \pi i \zeta_0}}{e^{2 \pi i \tau_j}-e^{2 \pi i \zeta_0}}
\right).  
\]
On the other hand, 
\[
\prod_{j=1}^4 (e^{2 \pi i \sigma_j}-e^{2 \pi i \zeta_0}) = 
(1-e^{2 \pi i \zeta_0}) \prod_{j=1}^3 (e^{2 \pi i \tau_j}-e^{2 \pi i \zeta_0})
\]
since $F_\zeta(\zeta_0) = 0$.
Let
\begin{equation}
a_2 \, u^2 + a_1 \, u + a_0
=
\left(\prod_{j=1}^4 (e^{2 \pi i \sigma_j}-u) - (1-u) \prod_{j=1}^3 (e^{2 \pi i \tau_j}-u)\right)/u.  
\label{eq:quad2}
\end{equation}
Then $e^{2 \pi i \zeta_0}$ is a solution of $a_2 \, u^2 + a_1 \, u + a_0=0$.  
Let $u_0 = e^{2 \pi i \zeta_0}$ and $u_0^\prime$ be another solution of \eqref{eq:quad2}.  
Then we have
\begin{align}
&\frac{\prod_{j=1}^4 (e^{2 \pi i \sigma_k}-u_0)}{u_0^2}\, F_{\zeta\zeta}(\zeta_0, \boldeta) 
\nonumber
\\
&= 
\frac{\pi\, i}{u_0}  \!  \left(\!
-\prod_{j=1}^3(e^{2 \pi i \tau_j}\!-u_0)+ 
\sum_{j=1}^4\prod_{1\leq k \leq 4, k\neq j}\!\!\!\!\!\!(e^{2 \pi i \sigma_k}\!-u_0)
-\!\!\!
\sum_{j=1}^3(1-u_0)\!\!\!\!\!\!\prod_{1 \leq k \leq 3, k \neq j}\!\!\!\!\!\!(e^{2 \pi i \tau_k}-u_0)\!\!
\right)\nonumber
\\
&=
\pi \, i \, \left(
3 \, a_2 \, u_0+ 
2\, a_1 + a_0 \, u_0^{-1}
\right) 
=
\pi \, i   \, \left(
a_2 \, u_0- 
a_0 \, u_0^{-1}
\right) \nonumber
=
a_2 \, \pi \, i   \, \left(
u_0-u_0^\prime
\right) 
\\
&
=
\pm \pi \, i \, \sqrt{a_1^2 - 4\, a_0 \, a_2}
=
\pm4 \, \pi \,  \sqrt{-\det G}.  
\label{eq:gram}
\end{align}
Here we use  
\[
a_1^2 - 4\, a_0 \, a_2 = 16 \, \det G  
\]
which is obtained by direct computation.  
\end{proof}
Recall  that the volume ${\rm Vol}(T)$ is given in \cite{MY} and \cite{U} as follows.  
\begin{multline}
{\rm Vol}(T)
=
\\
\left|2 \, \pi\, \left(F(\zeta_0, \boldeta) + \delta(\eta_a, \eta_b, \eta_e) 
+ \delta(\eta_a, \eta_c, \eta_f)
+ \delta(\eta_b, \eta_d, \eta_f)
+ \delta(\eta_c, \eta_d, \eta_e) 
\right)\right|.  
\label{eq:vol}
\end{multline}
Combining \eqref{eq:remain}, \eqref{eq:limit}, \eqref{eq:recont}, \eqref{eq:adelta}, \eqref{eq:ag}, \eqref{eq:gram} and \eqref{eq:vol}, 
we get
\begin{multline*}
\left|\left\{\begin{matrix}
a_r & b_r & 
e_r \\
d_r & c_r & f_r
\end{matrix}\right\}_{q = \xi_r^2}^{RW}
\right|
\underset{r\to\infty}{\sim}
\\
2\pi \, \tfrac{r^2}{16} \,
\tfrac{16}{r^3} \left|
\tfrac{e^{2 \pi i \zeta_0}}{\prod_{j=1}^4\sqrt{1-\frac{e^{2 \pi i \zeta_0}}{e^{2 \pi i \sigma_j}}}} 
\sqrt{\tfrac{2\, \pi}{r \left|F_{\zeta\zeta}(\zeta_0, \boldeta)\right|}} \, 
\exp\left(\frac{r}{2 \pi}  {\rm Vol}(T)
+ O(\tfrac{1}{r})
\right)\right|
\underset{r\to\infty}{\sim}
\\
\frac{\sqrt{2} \,  \pi}{r^{3/2}\, \sqrt[4]{-\det G}}
\,e^{\frac{{\rm Vol}(T) \, r}{2\, \pi} + O(\tfrac{1}{r})}.  
\end{multline*}
\section{Symmetry}
A widely applicable version of the volume conjecture was introduced by the first author and T. Yang in \cite{CY}.  
The original volume conjecture is for knots and links in $S^3$, while this new version works well for the Turaev-Viro invariants  \cite{TV} and the Witten-Reshetikhin-Turaev invariants  \cite{RT} of  3-manifolds.   
Here the Turaev-Viro invariant include the extension to cusped 3-manifolds and hyperbolic 3-manifolds with totally geodesic boundaries (see \cite{BP,CY}).  
Let $\xi_r = \exp(2\pi\sqrt{-1}/r)$ for a positive odd integer $r \geq3$.  
Substituting $q = \xi_r^2$ instead of $q=\xi_r$ and doing numerical experiments for some manifolds, they observe that
 the Turaev-Viro invariant and the Witten-Reshetikhin-Turaev invariant have exponential growth with respect to $r$ and the leading asymptotic is determined by the hyperbolic volume of the manifold.  
\begin{conjecture}[Chen-Yang \cite{CY}]
Let $M$ be a hyperbolic manifold and let ${\rm Vol}(M)$ be its hyperbolic volume.
\begin{enumerate}
\item[(1)]
Let $M$ be a 3-manifold which is closed, cusped, or with totally geodesic boundary, and
let $TV_r(M)$ be the Turaev-Viro invariant \cite{TV,CY,BP} with $q = \xi_r^2$ for odd $r$.  
Then
\[
\lim_{r\to \infty,\  \text{$r$\,:\,odd}}
\frac{2\, \pi}{r} \log\left|TV_r(M)\right|
=
{\rm Vol}(M).  
\]
\item[(2)]
Let $W\!RT_r(M)$ be the Witten-Reshetikhin-Turaev invariant \cite{W,RT,BHMV} with $q = \xi_r^2$ for odd $r$.  
Then
\[
\lim_{r\to \infty,\  \text{$r$\,:\,odd}}
\frac{4\, \pi\, \sqrt{-1}}{r} \log W\!RT_r(M)
=
{\rm CS}(M) + {\rm Vol}(M)\, \sqrt{-1}  \quad \mod \pi^2 \, {\mathbf Z}
\]
where ${\rm CS}(M)$ is the Chern-Simons invariant of $M$ multiplied by $2\pi^2$.  
\end{enumerate}
\label{conj:CY}
\end{conjecture}
\par
Conjecture \ref{conj:CY} is proved for some cases, one is for 3-manifolds obtained by integral surgeries along the figure-eight knot by T. Ohtsuki in \cite{O2}, and another one is for the complements of the figure-eight knot and the Borromean rings by R. Detcherry, E. Kalfagianni and T. Yang in \cite{DKY} and those for fundamental shadow links in \cite{BDKY}.  
Since the Turaev-Viro invariant is constructed by using the quantum $6j$ symbols,  this conjecture suggests that the terms of asymptotic expansion  of the  quantum $6j$ symbol are all expressed by some geometric data of the corresponding tetrahedra.  
\par

We are considering  a non-oriented tetrahedron.  
In this case, we can think of its volume, area of faces, lengths of edge and dihedral angles along edges.   
But there is no canonical way to assign the signs of these contents.  
So  formulas of geometric contents of a tetrahedron must have certain symmetry for changing the signs of the lengths and angles.   
For example, the entries of the Gram matrix is given in terms of cosines of the dihedral angles and they don't depend on the signs of dihedral angles.  
In the construction of the Turaev-Viro invariant, 
the normalized quantum $6j$ symbol
$\left|\begin{matrix}
a & b & e \\ d & c & f\end{matrix}\right|$ 
is assigned to each tetrahedron instead of 
$\left\{\begin{matrix}
a & b & e \\ d & c & f\end{matrix}\right\}^{RW}$ where 
$$
\left|
\begin{matrix}
a & b & e \\ d & c & f
\end{matrix}
\right|
=
\sqrt{-1}^{-2(a+b+c+d+e+f)} \, 
\left\{\begin{matrix}
a & b & e \\ d & c & f\end{matrix}\right\}^{RW}.  
$$
So it is natural to expect certain symmetry about each parameter of  $\left|\begin{matrix}
a & b & e \\ d & c & f\end{matrix}\right|$ reflecting the change of the sign of the corresponding dihedral angle.  
Here we consider such symmetry for the following continuous version of the quantum $6j$ symbol.  
Let
\begin{multline*}
F_r(\theta_a, \theta_b, \theta_c, \theta_d, \theta_e, \theta_f) 
= 
\\
2 \, (-1)^{\frac{r+1}{2}} \sin\frac{2\pi}{r} 
e^{\frac{2\pi i}{r}\left(d_3(a_r, b_r, e_r)+d_3(a_r, c_r, f_r) + d_3(b_r, d_r, f_r)+d_3(c_r, d_r, e_r)\right)}
\times
\\
%\times
%\\
{\widetilde\Delta}_r(a_r, b_r,e_r)\,
{\widetilde\Delta}_r(a_r, c_r, f_r)\,
{\widetilde\Delta}_r(b_r, d_r, f_r)\,
{\widetilde\Delta}_r(c_r, d_r,e_r)
\int_{-\infty}^\infty\psi_r(z)\, \widetilde \alpha_r(z)\, dz,    
\end{multline*}
where $a_r = \frac{r(\pi-\theta_a)}{4\pi}-\frac{1}{2}$, $\cdots$, $f_r = \frac{r(\pi-\theta_f)}{4\pi}-\frac{1}{2}$.  
For the functions $\widetilde\alpha_r$, $\widetilde\Delta_r$, $\psi_r$, refer \eqref{eq:contalpha}, \eqref{eq:contdelta} and \eqref{eq:psi} respedtively.  
According to the argument in Appendix, 
the asymptotics of 
$
\left|
\begin{matrix}
a_r & b_r & e_r \\ d_r & c_r & f_r
\end{matrix}
\right|
$
with respect to $r$
is equal to that of
$
F_r(\theta_a,\theta_b$, $\theta_c$, $\theta_d$, $\theta_e$, $\theta_f) 
$
where $a_r = \frac{r(\pi-\theta_a)}{4\pi} - \frac{1}{2}$, $\cdots$, $f_r = \frac{r(\pi-\theta_f)}{4\pi} - \frac{1}{2}$.  
\begin{conjecture}[Symmetry of asymptotics]
Let $T$ be a hyperbolic truncated tetrahedron and  $\theta_a$, $\theta_b$, $\theta_c$, $\theta_d$, $\theta_e$, $\theta_f$ be dihedral angles at edges $a$, $\cdots$, $f$ in Figure \ref{fig:tet}.
%\textcolor{red}{
%Let  
%\[
%\eta_a = \frac{1}{2}-\frac{\theta_a}{2\pi}, \quad
%\eta_b = \frac{1}{2}-\frac{\theta_b}{2\pi}, \quad
%\cdots,\quad 
%\eta_f = \frac{1}{2}-\frac{\theta_f}{2\pi},   
%\]
%and 
%let 
%\[
%\eta_a^\prime = \frac{1}{2}+\frac{\theta_a}{2\pi}.  
%\]
Then the asymptotic expansions of 
$F_r(\theta_a$, $\theta_b$, $\theta_c$, $\theta_d$, $\theta_e$, $\theta_f)$
and 
$F_r(-\theta_a$, $\theta_b$, $\theta_c$, $\theta_d$, $\theta_e$, $\theta_f)$
with respect to $r$ are  expected to agree at all orders of $\frac{1}{r}$.
\label{conj:symmetric}
\end{conjecture}
This conjecture explains the very mysterious big cancellation happening in the Chen-Yang's Volume Conjecture of the Turaev-Viro invariants.  
\par
\begin{remark}
If $a_r$, $b_r$, $c_r$, $d_r$, $e_r$, $f_r$ are admissible and $a_r^\prime$, $b_r^\prime$, $c_r$, $d_r$, $e_r$, $f_r^\prime$ are also admissible,  then we have the following identity proved in \cite{DKY}.   
$$
\left|
\begin{matrix}
a_r & b_r &e_r \\ 
d_r & c_r & f_r
\end{matrix}
\right|_{q=\xi_r^2}
=
\left|
\begin{matrix}
a_r^\prime & b_r^\prime &e_r \\ 
d_r & c_r & f_r^\prime
\end{matrix}
\right|_{q=\xi_r^2}
$$
where $a_r^\prime = r\,\frac{\pi+\theta_a}{4\pi} - \frac{1}{2}$, i.e. $\theta_a^\prime = - \theta_a$.  
Applying our Conjecture 3 three times predicts that those $6j$-symbols have the same expansion, but in this case the identity shows that they are completely equal.
\end{remark}
\begin{remark}
Let 
\[
F_r(\theta_a, \theta_b, \theta_c, \theta_d, \theta_e, \theta_f)
\underset{r\to \infty}\sim
(-1)^{\frac{r+1}{2}} \, 
A \,
e^{\frac{r}{2\,\pi}B} \,
\left( 1 + C_1 \, \frac{1}{r} + C_2 \, \frac{1}{r^2} + \cdots \right)
\]
be the asymptotic expansion of 
$F_r(\theta_a, \theta_b, \theta_c, \theta_d, \theta_e, \theta_f)$.  
Then we showed that 
\[
A= \frac{\sqrt{2}}{r^{3/2} \, \sqrt[4]{-\det G}}, \qquad
B = {\rm Vol}(T),
\]  
where ${\rm Vol}(T)$ is the volume of the tetrahedron $T$ with dihedral angles $\theta_a$, $\cdots$, $\theta_f$ and $\det G$ is the Gram matrix of $T$. 
The Conjecture \ref{conj:symmetric} means that the coefficients $A$, $B$, $C_1$, $C_2$, $\cdots$ are all symmetric with respect to the switching of the sign of  dihedral angles.  
The volume function ${\rm Vol}(T)$ and the Gram matrix $\det G$ are not changed by changing the sign of any dihedral angle, the above conjecture is true up to the second leading term.  
\end{remark}
We computed  the term $C_1$.    
%Let $T_\theta$ and $T_{-\theta}$ be the regular tetrahedra with dihedral angles $\theta$ and $-\theta$ respectively, then $T_{-\theta}$ is obtained by applying the flips in Conjecture~\ref{conj:symmetric} to every edges.  
%The tetrahedra $T_{\theta}$ and $T_{-\theta}$ are both ultra ideal for  $-\pi/3 < \theta < \pi/3$, and  
Let $C_1(\theta_a)$ denote the therm $C_1$  with 5 fixed parameters ($\theta_b=\frac{2}{25}\pi$, $\theta_c=\frac{3}{25}\pi$, $\theta_d=\frac{4}{25}\pi$, $\theta_e = \frac{1}{5}\pi$, $\theta_f = \frac{6}{25}\pi$).  
We computed $C_1(\theta_a)$ numerically as the following table and graph. 
%$C_1(\theta_a)$ diverges at $\theta_a = \pm\frac{16}{25}\pi$, $\pm\frac{18}{25}\pi$ where a vertex of the corresponding tetrahedron is ideal.  
It suggests that $C_1(\theta_a) = C_1(-\theta_a)$.  
\begin{table}[htb]
\begin{center}
\begin{tabular}{|c|ccccccc|}
\hline
$\theta_a$\vphantom{$\frac{2^{1}}{1_{1_1}}$} & $-\frac{3}{5}\pi$ &$-\frac{2}{5}\pi$ & $-\frac{1}{5}\pi$ & $0$ & $\frac{1}{5}\pi$ & $\frac{2}{5}\pi$ & $\frac{3}{5}\pi$
\\
\hline
$C_1(\theta_a)$ \vphantom{$\frac{2^{1}}{1_{1_1}}$}& 
$4.85928$ & $2.33391$ & $1.98659$ & $1.90677$ & 
$1.98659$ & $2.33391$ &  $4.85928$
\\
\hline
\end{tabular}
\\{\ }\\
\end{center}
\caption{Values of \,$C_1(\theta_a)$. }
\end{table}
\begin{figure}[htb]
\begin{center}
\qquad\qquad\small$C_1(\theta_a)$
\\
\epsfig{file=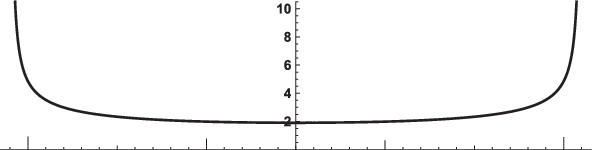, scale=0.9}
\raisebox{1mm}{\small$\theta_a$}
\\
\scriptsize{
$-\frac{3}{25}\pi$ \qquad \qquad \qquad\qquad \ \ $-\frac{1}{25}\pi$ \quad\qquad\quad\  $0$ \ \ \qquad\qquad $\frac{1}{25}\pi$    \quad\qquad\qquad \qquad\qquad\ $\frac{3}{25}\pi$\ \ \ \ \ \ \ }
\end{center}
\caption{Graph of $C_1(\eta_a)$.}
\end{figure}
\section*{Appendix. Comparison of \,$\widehat h_r(0)$ and the sum of $\widetilde\alpha_r(z)$}
Let $\beta_r(x) = \frac{1}{r} \, \log\widetilde\alpha_r(\frac{r\,x-3}{2})$.  
Then $\beta_r(x)$  is a real analytic function defined on the interval $(u, v)$ with $u = \frac{2m_r+3}{r}$ and $v=\frac{2M_r+3}{r}$ where $m_r$ and $M_r$ are the lower and upper bounds in the line below \eqref{eq:alpha}.  
Note that $u$ and $v$ do not depend on $r$.  
Such function $\beta_r(x)$ exists since $\widetilde\alpha_r(\frac{r\,x-3}{2})$ is a  positive real analytic function.  
Thanks to Lemma 2,  $\beta_r(x)$ has a limiting function $\beta(x)=\lim_{r\to\infty}\beta_r(x)$,   
which is a real analytic function which takes its maximum at $x = x_0 \in (u, v)$.  
Recall that
$
h_r(z) = \psi_r(\frac{2z+3}{r}) \, \widetilde \alpha_r(z)  
$.
Then, by applying the Poisson summation formula, we have
\[
\sum_{z \in \mathbf Z} h_r(z) =
\sum_{m \in \mathbf Z} \widehat h_r(m), 
\]
and 
\[
\begin{aligned}
\widehat h_r(m)
%=
%\int_{\mathbf R}e^{\pi i r m x} \,  f(x) \, dx
&=
\int_{\mathbf R} e^{-2 \pi i m z}\, \psi\left(\frac{2z+3}{r}\right)\, \widetilde\alpha_r(z) \, dz
%\\&
=
\int_{\mathbf R} e^{-2 \pi i m \frac{r x-3}{2}}\, \psi(x)\, \widetilde\alpha_r\left(\frac{r\,x-3}{2}\right) \, \frac{r}{2} \, dx
\\
&=
(-1)^m \, \frac{r}{2} \, \int_{\mathbf R} \psi(x) \, e^{r(\beta_r(x)- \pi i m x)} \, dx.
\end{aligned}
\]
%\par
%
\begin{lemma}
There are positive constants $c''$ and $\delta$ satisfying
\begin{equation}
\frac{\left|\widehat h_r(1)\right|}{\left|\widehat h_r(0)\right|} < c'' \, e^{-r\, \delta}.  
\label{eq:limit2}
\end{equation}for some $\delta>0$.  
\label{lem:large}
\end{lemma}
%
%Now we prove Lemma \ref{lem:large}. 
%
\begin{proof}%[Proof of Lemma \ref{lem:large}]
At first, we extend $\beta(x)$ to a complex-valued function around $x = x_0$.  
Let $\beta_{r, 1}(x)$ $= \beta_r(x) - \pi \,i \, x$.  
Choose $\varepsilon' > 0$ so that 
$\beta(x)-\pi\, i\, x $ doesn't have a critical point in the $\varepsilon'$ neighborhood $B(x_0, \varepsilon')$ of $x_0$ in 
$\mathbf C$ and 
$[u+2\varepsilon, v-2\varepsilon] \cap B(x_0, \varepsilon') 
= (x_0-\varepsilon', x_0 + \varepsilon')$ 
where $\varepsilon$ is the small positive number introduced in the definition of $\Psi$ in \eqref{eq:psi}.  
Let $x_0^{(r)}$ be the maximum point of $\beta_r(x)$, then $x_0^{(r)}$ is close to $x_0$ when $r$ is sufficiently large and 
$$
\beta_r(x) = \beta_r(x_0^{(r)}) - a_r \, (x-x_0^{(r)})^2 + O((x-x_0^{(r)})^3)
$$
in a small neighborhood of $x_0^{(r)}$.  
For $x\in {\mathbf C}$ near $x_0^{(r)}$, the contours of ${\rm Re} \, \beta_r(x)$ and ${\rm Re}\, \beta_{r,1}(x)$ are given in the following figure.  
Choose  
$x_1 \in (x_0-\varepsilon', x_0)$ 
so that the contour 
$C_{r,1}$ 
corresponding to the value 
${\rm Re}\,\beta_{r, 1}(x_1)$ 
connects $x_1$ to a real point $x_1^{(r)}\in (x_0, x_0 + \varepsilon)$  in $B(x_0, \varepsilon')$ and for any sufficiently large $r$ as in the figure.  
Please note that ${\rm Re} \, \beta_{r,1}(x) = {\rm Re} \, \beta_r(x)$ for $x \in {\mathbf R}$. % and so $x_1^{(r)}$ doesn't depend on $m$.  
%Such $x_1$ exists since $\beta_r(x)$ converges to $\beta(x)$ and   
%the contour $C_{r,m}$ $(|m| \geq 2)$ is in the region enclosed by $C_{r, 1}$ %and  $C_{r, -1}$.  
Then, due to $\beta_r(x_1) \to \beta(x_1)$, $\beta_r(x_0^{(r)}) \to \beta(x_0)$ when $r\to \infty$ and $\beta(x_1) \neq \beta(x_0)$,  there is $\delta > 0$ such that 
 $\beta_r(x_1) < \beta_r(x_0^{(r)})-\delta$ for any sufficiently large $r$.
\begin{figure}[htb]
$$
\begin{matrix}
\raisebox{3.5cm}{$\mathbf C$}
\hspace{2.1cm}\raisebox{1.56cm}{\small$x_0^{(r)}$}\hskip-2.6cm
\epsfig{file=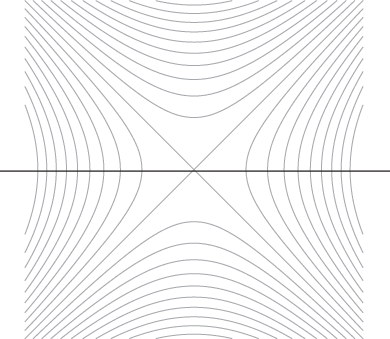, scale=0.7}\raisebox{1.8cm}{$\mathbf R$}&\qquad 
\raisebox{3.5cm}{$\mathbf C$}
\hspace{1.6cm}\raisebox{2.2cm}{\small $x_1$ \ \ \ $x_0^{(r)}$ \ \ $x_1^{(r)}$}
\hskip-1.4cm \raisebox{7mm}{$B(x_0, \varepsilon)$}
\hskip-0.5cm\raisebox{1.5cm}{$C_{r,1}$}\hskip-3.2cm
\epsfig{file=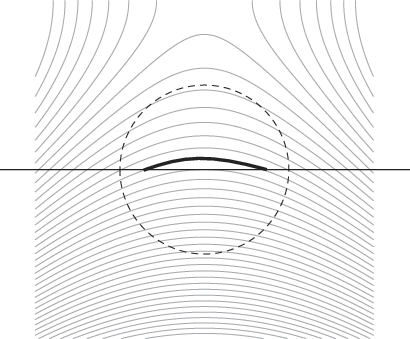, scale=0.7}\raisebox{1.8cm}{$\mathbf R$}
\\
\text{Contours of ${\rm Re}\,\beta_r(x)$}
 \quad&\quad 
 \text{Contours of ${\rm Re}\,\beta_{r,1}(x)$}
\end{matrix}
$$
\end{figure}
\par
Now recall that 
\begin{align*}
\widehat h_r(0) 
&=
%\frac{r}{2} \, \int_{\mathbf R} h_r(x) \, dx
%=
%\frac{r}{2} \, \int_{u+\varepsilon}^{v-\varepsilon} \psi(x) \, \widetilde\alpha_r\left(\frac{r\,x-3}{2}\right) \, dx
%\\
%&=
\frac{r}{2} \, \int_{u+\varepsilon}^{v-\varepsilon} \, \psi(x) \, e^{r \beta_r(x)} \, dx \underset{r\to\infty}{\sim}
\text{const.} \, r^{-5/2}\, e^{r\, \beta(x_0)} ,
\\
\left|\widehat h_r(1)\right|
&=
% \frac{- r}{2} \, \int_{u+\varepsilon}^{v-\varepsilon} h_r(x) \, dx 
% =
\frac{r}{2} \, \left|\int_{u+\varepsilon}^{v-\varepsilon} \psi(x) \, e^{r (\beta_r(x)-\pi i x)} \, dx\right|
%\\&
\leq
\frac{r}{2} \, \int_{C'}\psi(x) \, e^{r (\operatorname{Re}(\beta_r(x)-\pi i x))} \, dx
\\&
\leq
\frac{r}{2} \, \int_{C'}\psi(x) \, e^{r (\operatorname{Re}\,\beta_r(x_1))} \, dx
=
\frac{r}{2} \, \ell(C') \, e^{r\, \operatorname{Re}\, \beta_r(x_1)}
\end{align*}
where $C' = [u, x_1]\cup C_{r, 1} \cup [x_1^{(r)}, v]$ and $\ell(C')$ is the length of $C'$.  
Since  $\operatorname{Re}\, \beta(x_1) < \operatorname{Re}\,\beta(x_0)-\delta$, the relation \eqref{eq:limit} holds for some positive constant $c''$.  
\end{proof}
%
%\begin{lemma}
%There is a constant $C$ not depending on $r$ such that 
%\begin{equation}
%\widehat h_r(0) > C\, \sqrt{r} \, e^{r \, \beta_r(x_1)}
%\label{eq:C}
%\end{equation}
%for large $r$.  
%\end{lemma}
%%
%\begin{proof}
%Since 
%$\widehat h_r(0) \sim \frac{r}{2} \, \frac{e^{r \beta_r(x_0^{(r)})}}{\sqrt{-r \, \beta_r''(x_0^{(r)})}}$, there is a constant $C_0$ such that 
%$
%\frac{2 \, \widehat h_r(0)\,\sqrt{-\beta_r''(x_0^{(r)})}}{\sqrt{r} \, e^{r \beta_r(x_0^{(r)})}}$ $>C_0$
%for any large $r$.   
%Put $C = \frac{C_0}{2 \, \sqrt{-\beta''_r(x_0^{(r)})}}$, then we proved \eqref{eq:C}.  
%\end{proof}
%
\begin{lemma}
For $m>1$,  $\widehat h(m)$ satisfies
\begin{equation}
\left|\widehat h_r(m)\right| \leq \frac{D}{m^2}\, \left|\widehat h_r(1)\right| .
\label{eq:c}
\end{equation}
for some constant $D$ not depend on $m$. 
\end{lemma}
\begin{proof}
We have $\operatorname{Re}\,(\beta(z) - \pi\, m\, i\, z) = \operatorname{Re}\,(\beta(z) - \pi\,  i\, z)$ on $[u, v]$ and $\operatorname{Re}\,(\beta(z) - \pi\, m\, i\, z) \leq \operatorname{Re}\,(\beta(z) - \pi\,  i\, z)$ on $C_{r, 1}$.  
Now we reformulate $\widehat h_r(m)$ as
\begin{align*}
&\widehat h_r(m)
=
(-1)^m \, \frac{r}{2} \, \int_{C'}e^{-\pi i r m }\,  \psi(x) \, e^{r \beta_r(x)x} \, dx
\\
&=
\frac{(-1)^{m+1} \,r}{2\,\pi \, i\, m} \, \int_{C'}e^{-\pi i r m }\,  
\left(\psi'(x) + r\, \beta'_r(x)\right)\, e^{r \beta_r(x)x} \, dx
\\
&=
\frac{(-1)^{m+1} \,r}{4\,\pi^2 \,  m^2} \, \int_{C'}e^{-\pi i r m }\,  
\left(\psi''(x) + r\, \beta''_r(x)
+r\, \left(\psi'(x) + r\, \beta'_r(x)\right)\, \eta'_r(x)\right)\, e^{r \beta_r(x)x} \, dx
\end{align*}
then there is a constant $d$ satisfying
\[
\left|
e^{-\pi i r m }\,  
\left(\psi''(x) + r\, \beta''_r(x)
+r\, \left(\psi'(x) + r\, \beta'_r(x)\right)\, \eta'_r(x)\right)
\right| < d.
\]
Hence $\left| \widehat h_r(m)\right| < 
\frac{d\, r}{4 \pi^2 m^2} \, \left|\widehat h_r(1)\right|$.  
%
%
%\begin{align*}
%\left|\widehat h_r(m)\right|
%&=
%\left| \frac{r}{2} \, \int_{\overline C_{r.,m}}e^{\pi i r m x} \, \psi_r(x) \, e^{r \beta_r(x)} \, dx\right|
%%\\&
%\leq
%\frac{r}{2} \, \int_{\overline C_{r.,m}}\left| e^{\pi i r m x} \, \psi_r(x) \, e^{r \beta_r(x)} \right|\, dx
%\\
%&=
%\frac{r}{2} \, \int_{\overline C_{r.,m}} \psi_r(x) \, e^{r \,{\rm Re}\,\beta_r(x)} \, dx
%%\\&
%\leq
%\frac{r}{2} \, \int_{\overline C_{r.,m}} \psi_r(x) \, e^{r \beta_r(x_1)} \, dx
%%\\&
%\leq
%\frac{r}{2} \,L_r \,  e^{r \beta_r(x_1)}
%\end{align*}
%where $L_r$ is the length of $\overline C_{r,m}$ and $\psi_r(x)=1$ for $x \in C_{r,m}$.  
%So we proved \eqref{eq:c} by setting $c = \frac{\underset{r}\sup\, L_r}{2}$.    
%\medskip
\end{proof}
\begin{proof}[Proof of Lemma \ref{lem:remain2}]
The above properties also holds for $\widehat h_r(-1)$ and $\widehat h_r(m)$ for negative $m$.  
So we have $\left|\sum_{m\neq0} \widehat h_r(m)\right| \leq M \,
\left|\widehat h_r(1)\right|$ for some constant $M$  since $\sum_{m \neq 1}\frac{1}{m^2}$ converges. 
Hence Lemma \ref{lem:large} implies \eqref{eq:limit}.  
\end{proof}
%
%\section*{Appendix B. Comparison of $h(z)$ and $\widetilde \alpha(z)$}
%
%In this appendix, we show the following.
%
Now we prove Lemma \ref{lem:remain}.  
\begin{proof}[Proof of Lemma \ref{lem:remain}]
Let $x_1$, $x_1^{(r)}$ be points on $[u,v]$ given in the proof of Lemma \ref{lem:large}.  
Then $\Psi(x) = 1$ for $x \in [x_1, x_1^{(r)}]$.  
Therefore $\widetilde \alpha_r(x) - h_r(x)=0$ on $[x_1, x_1^{(r)}]$.  
Since 
$|h_r(x)| \leq |\widetilde\alpha_r(x)|$ and $|\widetilde\alpha(x)|$ is increasing for $x \leq x_1$ and decreasing for $x \geq x_1^{(r)}$.  
Hence 
\[
|A_r - B_r| \leq 2\, r \,  |\widetilde\alpha_r(x_1)| = 2\, r\, e^{r
\, \beta_r(x_1)} < 2\, r\, e^{r\,( \beta_r(x_0^{(r)})-\delta)}.
\]  
By putting $\delta'=\delta$, 
this implies \eqref{eq:remain}.  
\end{proof}
\section*{Acknowledgements}
The first author would like to thank Francis Bonahon and Nicolai Reshetikhin for their strong interests, valuable discussions, as well as cordial guidance during past years. He also thank Giovanni Felder for his large interests, strong support and enlightening discussions, and Rinat Kashaev, Tomotada Ohtsuki, Shing-Tung Yau for their strong interests, hospitality, encouragement and inspiring discussion. 
Finally he would like to thank Tian Yang for his proof reading, discussion and many valuable suggestions.
%\end{acknowledgements}

\begin{small}
%\end{document}
\address{
   Qingtao Chen, 
   Division of Science,
   New York University Abu Dhabi, PO Box 129188
   Abu Dhabi, United Arab Emirates,} 
   \email{chenqtao@nyu.edu}
   \\
% Important: Do not put any empty line here.

\address{
Jun Murakami, 
   Department of Mathematics,
Faculty of Science and Engineering,
Waseda University, 
3-4-1 Ohkubo, Shinjuku-ku, Tokyo 169-8555
Japan,}
\email{murakami@waseda.jp}
\end{small}
\end{document}